\documentclass[12pt]{article}
\usepackage{amssymb, amsmath, latexsym, amsbsy}
\usepackage{amsfonts}
\usepackage{dsfont}
\usepackage[T1]{fontenc}
\usepackage[all]{xy}
\usepackage{amscd}
\usepackage{graphicx}
\usepackage{latexsym}
\usepackage{pxfonts}
\usepackage{enumerate}
\usepackage{float}
\usepackage{color}
\usepackage[colorlinks]{hyperref}
\usepackage{fancybox}
\usepackage{mathtools, booktabs, tabularx, adjustbox}
\usepackage{fancyhdr}
\usepackage{hyperref}

\newtheorem{Definition}{Definition}
\newtheorem{Theorem}{Theorem}
\newtheorem{Lemma}{Lemma}
\newtheorem{Corollary}{Corollary}
\newtheorem{Remark}{Remark}
\newtheorem{Proposition}{Proposition}

\newenvironment{proof}[1][Proof]{\textbf{#1:} }{\ \rule{0.5em}{0.5em}}
\title{Geometry of left-invariant vector fields on Lie groups.}
\begin{document}
\maketitle	

	\begin{center}
		\author{ \textbf{M. L. Foka}$^{1}$, \quad \textbf{   R.P. Nimpa}$^{2}$,\quad \textbf{M. B. N.  Djiadeu}$^{3}$,\\
			\small{e-mail: $\textbf{1}.$ lanndrymarius@gmail.com, \quad
				$\textbf{2}.$ romain.nimpa@facsciences-uy1.cm, \\
				$\textbf{3}.$  michel.djiadeu@facsciences-uy1.cm,   \\
				University of Yaounde 1, Faculty of Science, Department of
				Mathematics,  P.O. Box 812, Yaounde, Republic of Cameroon.}}
	\end{center}
		
		\begin{abstract}
			
			We investigate the geometry of left-invariant vector fields on simply connected nilpotent Lie groups equipped with left-invariant Riemannian metrics. Exploiting the canonical identification between the Lie algebra $\mathfrak{g}$ and the space of left-invariant vector fields, we establish complete algebraic characterizations for several fundamental classes of vector fields, including Killing, one-harmonic, harmonic, conformal, and concurrent fields.
			
			Our main results reveal a striking rigidity phenomenon: on any nilpotent Lie group, the spaces of Killing, one-harmonic, and conformal vector fields all coincide precisely with the center of the Lie algebra. Moreover, we prove that no nontrivial concurrent vector fields exist in this setting. In contrast, harmonic vector fields form a proper subspace of the Lie algebra, characterized as those central vectors orthogonal to the derived algebra.
			

		\end{abstract}

	\textbf{keywords}: Lie group; Killing vector field; one-harmonic vector field; conformal vector field; concurrent vector field.

	
	\section{Introduction}
	
	The geometry of Lie groups endowed with left-invariant Riemannian metrics relies fundamentally on the canonical isomorphism between the Lie algebra and the space of left-invariant vector fields. This identification permits the explicit expression of the Levi-Civita connection, divergence, vector Laplacian, and Ricci tensor in terms of Lie brackets and some associated algebraic operators, following classical approaches by Milnor \cite{Milnor} and O'Neill \cite{ONeill}.
	
	The study of invariant harmonic vector fields trace back to González-Dávila and Vanhecke \cite{GonzalezDavila2002}, who fully characterized left-invariant unit harmonic vector fields on various Lie groups, including unimodular, Heisenberg, and Damek-Ricci types. Extensions to oscillator groups were given by Xu and Tan \cite{XuTan2019}. Harmonic and minimal vector fields have traditionally been examined by equipping the tangent bundle $TM$ with the Sasaki Riemannian metric, a natural construction derived from the base manifold $(M,g)$ that lifts its Riemannian structure (see \cite{ab, gcal1, gcal2}). This approach has been widely explored in the literature. However, when analyzing harmonicity conditions, it is equally pertinent to consider the complete lift metric $g^c$ on $TM$, which differs from the Sasaki lift by having a neutral signature $(n,n)$ \cite{be,ga,on,va}. This alternative framework has proven useful in various geometric contexts, including the study of harmonic sections and one-harmonic vector fields. Following this line of thought, Calviño-Louzao et al. \cite{cal} proceed to identify all left-invariant one-harmonic vector fields on three-dimensional Lie groups.
	
	In the context of left-invariant conformal vector fields, Cintra, Chen, and Neto \cite{Cintra2016} proved that on pseudo-Riemannian unimodular Lie groups, every conformal field is necessarily Killing. They also provided necessary conditions for conformal fields that are non-Killing in non-unimodular (notably Lorentzian) settings. More recently, Herrera and Origlia \cite{Herrera2023} provided a survey on invariant conformal Killing $p$-forms on left-invariant Lie groups, with emphasis on nilpotent cases.
	
	Our present work contributes to this tradition by establishing general results valid for nilpotent Lie groups of arbitrary dimension. More intuitively, the following
	Theorems are proved.
	
	
	
	\begin{Theorem}[Killing vector fields]\label{thm:main1}
		Let $G$ be a simply connected nilpotent Lie group of dimension $n$ equipped with a left-invariant Riemannian metric. Then the space of left-invariant Killing vector fields on $G$ coincides with the center of its Lie algebra $\mathfrak{g}$.
	\end{Theorem}
	
	\begin{Theorem}[One-harmonic vector fields]\label{thm:main2}
		Let $G$ be a simply connected nilpotent Lie group of dimension $n$ equipped with a left-invariant Riemannian metric. A left-invariant vector field on $G$ is one-harmonic if and only if it is a Killing vector field.
	\end{Theorem}
	
	
	\begin{Theorem}[Conformal vector fields]\label{thm:main3}
		Let $G$ be a simply connected nilpotent Lie group of dimension $n$ equipped with a left-invariant Riemannian metric.
		The space of left‑invariant conformal vector fields on $G$ coincides with the space of left‑invariant Killing vector fields.
	\end{Theorem}
	
	\begin{Theorem}[Concurrent vector fields]\label{thm:main4}
		Let $G$ be a simply connected Lie group of dimension $n$ equipped with a left-invariant Riemannian metric. There exists no nontrivial left-invariant concurrent vector field on $G$.
	\end{Theorem}
	\begin{Theorem}[Harmonic vector fields]\label{thm:harmonic-rigorous}
		Let $(G, h)$ be a simply connected nilpotent Lie group endowed with a left-invariant Riemannian metric, and let $\xi$ be a left-invariant vector field on $G$. Then $\xi$ is harmonic (i.e., $\Delta\xi = 0$) if and only if:
		\begin{enumerate}[(i)]
			\item $\xi$ is in the center of $\mathfrak{g}$ and ,
			\item $\xi \perp [\mathfrak{g}, \mathfrak{g}]$ (i.e., $\langle \xi, [u, v] \rangle = 0$ for all $u, v \in \mathfrak{g}$).
		\end{enumerate}
	\end{Theorem}
	
	The remainder of this manuscript is structured as follows. In Section 2, we recall foundational definitions and establish preliminary results. In Section 3, we prove our main theorems for nilpotent Lie groups of arbitrary dimension.
	\section{Preliminaries}
	


	Let $(M,g)$ be an $n$-dimensional Riemannian manifold and let $\nabla$ be its Levi-Civita connection.
	
	\begin{Definition}
		A vector field $X$ on $(M, g)$ is said to be a \emph{Killing field} if and only if
		\[
		\mathcal{L}_X g = 0,
		\]
		where $\mathcal{L}_X$ denotes the Lie derivative with respect to $X$.
	\end{Definition}
	
	\begin{Definition}\cite{Dod}
		The \emph{tension field} of a smooth map $\varphi: (M,g) \to (N,h)$ between Riemannian manifolds is defined by
		\[
		\tau(\varphi) = \operatorname{tr}(\nabla d\varphi),
		\]
		where $d\varphi$ is the differential of $\varphi$. The map $\varphi$ is called \emph{harmonic} if $\tau(\varphi) = 0$.
		
		A vector field $X$ on $(M, g)$ is called \emph{one-harmonic} if the local one-parameter group of diffeomorphisms $\{\varphi_t\}_{t\in I}$ it generates satisfies:
		\[
		\left.\frac{d}{dt}\,\tau(\varphi_t)\right|_{t=0}=0.
		\]
	\end{Definition}
	
	The following characterization theorem holds:
	
	\begin{Theorem}\cite{cal}\label{thm:one-harmonic-char}
		Let $X$ be a smooth vector field on $(M,g)$. The following conditions are equivalent:
		\begin{enumerate}
			\item $X$ is a one-harmonic vector field.
			\item $\operatorname{tr}(\mathcal{L}_X\nabla) = 0$.
			\item $X : (M,g)\to (TM,g^{c})$ is a harmonic map, where $g^c$ denotes the complete lift metric.
			\item $X$ is a Jacobi vector field along the identity map.
			\item $\Delta X = \operatorname{Ric}(X)$, where $\Delta$ is the vector Laplacian and $\operatorname{Ric}$ is the Ricci operator.
		\end{enumerate}
	\end{Theorem}

	\begin{Definition}\cite{gcal1}\label{def:harmonic}
		A vector field $X$ on $M$ is called \emph{harmonic} if
		\[
		\Delta X = 0,
		\]
		where $\Delta$ is the Hodge-de Rham Laplacian on vector fields, defined by:
		\[
		\Delta X = -\operatorname{tr}(\nabla^2 X) = -\sum_{i=1}^n (\nabla_{e_i}\nabla_{e_i} X - \nabla_{\nabla_{e_i}e_i} X),
		\]
		for any orthonormal local frame $\{e_1, \ldots, e_n\}$ on $(M, g)$.
	\end{Definition}
	
	\begin{Remark} 
		From Theorem \ref{thm:one-harmonic-char}, a vector field $X$ is called \emph{one-harmonic} if
		\[
		\Delta X = \operatorname{Ric}(X),
		\]
		where $\operatorname{Ric}$ is the Ricci operator. Harmonic vector fields correspond to the special case where the Ricci curvature vanishes in the direction of $X$. On Ricci-flat manifolds, the notions of harmonic and one-harmonic coincide.
	\end{Remark}

	\begin{Definition}\cite{Cintra2016}
		A vector field $X$ on $(M, g)$ is called a \emph{conformal vector field} if the local one-parameter group of diffeomorphisms $\{\phi_t\}$ generated by $X$ consists of local conformal transformations; that is, each local flow map $\phi_t$ satisfies
		\[
		\phi_t^* g = e^{2\psi_t} g
		\]
		for some smooth function $\psi_t \colon M \to \mathbb{R}$.
	\end{Definition}
	
	\begin{Proposition}\cite{w}
		The following statements are equivalent for a vector field $X$ on a Riemannian manifold $(M^n, g)$:
		\begin{enumerate}[(i)]
			\item $X$ is a conformal vector field,
			\item $\mathcal{L}_X g = 2h g$ for some function $h \in \mathcal{C}^{\infty}(M)$,
			\item $\mathcal{L}_X g = \dfrac{2\,\mathrm{div}(X)}{n} g$.
		\end{enumerate}
	\end{Proposition}
	
	\begin{Definition}
		A vector field $X$ on $(M, g)$ is called a \emph{concurrent vector field} if, for every vector field $Y$ on $M$, the following condition holds:
		\[
		\nabla_Y X = Y.
		\]
	\end{Definition}
	
	
	
	
	
	\section{Proofs of the main results}
	
	In this section, we prove our main results for nilpotent Lie groups of arbitrary dimension.
	We begin by recalling a fundamental fact from linear algebra.
	
	\begin{Lemma}\label{lem:nilpotent-skew}
		If $A \in M_n(\mathbb{R})$ is both skew-symmetric and nilpotent, then $A = 0$.
	\end{Lemma}
	
	\begin{proof}
		The Frobenius norm of $A$ is given by
		\[
		\|A\|_{\text{fr}}^2 := \operatorname{tr}(A^T A) = -\operatorname{tr}(A^2) = 0,
		\]
		where in the last equality we use the fact that $A$ is nilpotent (the trace of any power of a nilpotent matrix is zero). This implies that $A = 0$.
	\end{proof}
	

	
	

	Let $(G,h)$ be a Riemannian Lie group of dimension $n$, let $\mathfrak{X}_{\mathrm{inv}}(G)$ denote the set of left-invariant vector fields on $G$, let $(\mathfrak{g}, [\cdot ,\cdot], \langle\cdot,\cdot\rangle)$ be its Lie algebra, and let $\nabla$ be the Levi-Civita connection associated with $h$. We identify $\mathfrak{g}$ with $\mathfrak{X}_{\mathrm{inv}}(G)$ via the canonical isomorphism.
	
	Lemma \ref{lem:nilpotent-skew} will be applied in the context of the nilpotent Lie group $G$. Recall that in this setting, for any $X \in \mathfrak{g}$, the linear map $\operatorname{ad}_X : \mathfrak{g} \to \mathfrak{g}$ is nilpotent.\\

	For any left-invariant vector field $\xi$ on $G$, we consider the endomorphisms $\operatorname{L}_{\xi},\;\mathrm{R}_\xi,\; \operatorname{J}_\xi:\mathfrak{g}\longrightarrow \mathfrak{g}$ defined by:
	\[
	\operatorname{L}_{\xi}v = \nabla_{\xi} v,\quad \mathrm{R}_{\xi}v=\nabla_v\xi,\quad \operatorname{J}_{\xi}v=\operatorname{ad}^* _v\xi,
	\]
	where $\operatorname{ad}^*_v$ is the adjoint of $\operatorname{ad}_v$ with respect to the inner product $\langle\cdot,\cdot\rangle$.
	
	From the Koszul formula, we have the following identities:
	\[
	\operatorname{L}_\xi = \tfrac{1}{2} \left(\operatorname{ad}_\xi - \operatorname{ad}_\xi^* \right) - \tfrac{1}{2} \operatorname{J}_\xi,\quad
	\mathrm{R}_\xi = -\tfrac{1}{2} \left( \operatorname{ad}_\xi + \operatorname{ad}_\xi^* \right) - \tfrac{1}{2} \operatorname{J}_\xi.
	\]

	\subsection{Killing vector fields}
	The following lemma provides an algebraic characterisation of Killing vector fields on $G$ and is a preliminary step towards the proof of Theorem \ref{thm:main1}.
	
	\begin{Lemma}\label{lem:killing-char}
		A left-invariant vector field $\xi$ on $(G,h)$ is a Killing field if and only if $\operatorname{ad}_\xi$ is skew-adjoint, i.e., 
		\[
		\operatorname{ad}_\xi=-\operatorname{ad}_\xi^*.
		\]
	\end{Lemma}
	
	\begin{proof}
		For $u,v \in \mathfrak{g}$,
		\[
		(\mathcal{L}_\xi h)(u,v) = \xi h(u,v) - h([\xi,u],v) - h(u,[\xi,v]).
		\]
		By left-invariance of $h$, the function $h(u,v)$ is constant, so $\xi h(u,v) = 0$. Thus,
		\[
		(\mathcal{L}_\xi h)(u,v) = - h([\xi,u],v) - h(u,[\xi,v]).
		\]
		The condition $\mathcal{L}_\xi h = 0$ is equivalent to:
		\[
		\langle\operatorname{ad}_\xi (u),v\rangle + \langle u,\operatorname{ad}_\xi( v)\rangle = 0 \quad \forall u,v \in\mathfrak{g}.
		\]
		Using the adjoint definition, this is equivalent to $\operatorname{ad}_\xi = -\operatorname{ad}_\xi^*$.
	\end{proof}
	
	
	\begin{proof}[Proof of Theorem \ref{thm:main1}]
		Let $\xi \in \mathfrak{g}$. By Lemma \ref{lem:killing-char}, $\xi$ defines a left-invariant Killing vector field if and only if
		\[
		\operatorname{ad}_\xi = -\operatorname{ad}_\xi^*.
		\]
		Set $A := \operatorname{ad}_\xi : \mathfrak{g} \to \mathfrak{g}$. Since $\mathfrak{g}$ is nilpotent, the endomorphism $A$ is nilpotent. Being both nilpotent and skew-adjoint, we conclude by Lemma \ref{lem:nilpotent-skew} that $A = 0$, so $\xi$ lies in the center of $\mathfrak{g}$.
		
		Conversely, if $\xi$ belongs to the center of $\mathfrak{g}$, then $\operatorname{ad}_\xi = 0$, which is trivially skew-adjoint. Thus, $\xi$ defines a left-invariant Killing vector field.
	\end{proof}
	
	\subsection{One-harmonic vector fields}
	As preparation for the proof of Theorem \ref{thm:main2}, the subsequent lemma offers an algebraic viewpoint of one-harmonic vector fields on $G$
	
	\begin{Lemma}\label{lem:one-harmonic-char}
		A left-invariant vector field $\xi$ on $(G,h)$ is one-harmonic if and only if
		\begin{equation}\label{eq:one-harmonic}
			\sum_{i=1}^n(\operatorname{ad}_{v_i}^*+\operatorname{J}_{v_i})\circ \operatorname{ad}_\xi(v_i)- \tfrac{1}{2}\operatorname{ad}_{\xi}\circ \operatorname{ad}_{v_i}^*(v_i)=0,
		\end{equation}
		where $\{v_1,\cdots,v_n\}$ is an orthonormal basis of $\mathfrak{g}$ with respect to the inner product $\langle\cdot,\cdot\rangle.$
	\end{Lemma}
	
	\begin{proof}
		For all $i\in \{1,\cdots, n\},$ we have:
		\begin{align*}
			(\mathcal{L}_\xi \nabla)(v_i, v_i)
			&= [\xi, \operatorname{L}_{v_i} v_i] - \operatorname{L}_{[\xi, v_i]} v_i - \operatorname{L}_{v_i} [\xi, v_i]\\
			&= \operatorname{ad}_\xi \circ \operatorname{L}_{v_i}(v_i) - \mathrm{R}_{v_i}([\xi, v_i]) - \operatorname{L}_{v_i}\circ \operatorname{ad}_\xi(v_i)\\
			&= \operatorname{ad}_\xi \circ \Bigl(\tfrac{1}{2}(\operatorname{ad}_{v_i} - \operatorname{ad}_{v_i}^*) - \tfrac{1}{2} \operatorname{J}_{v_i}\Bigr)(v_i)
			-\Bigl(\operatorname{R}_{v_i} + \operatorname{L}_{v_i}\Bigr)\circ \operatorname{ad}_\xi(v_i) \\
			&= ( \operatorname{ad}_{v_i}^* + \operatorname{J}_{v_i} ) \circ \operatorname{ad}_\xi(v_i)- \tfrac{1}{2}\operatorname{ad}_{\xi}\circ \operatorname{ad}_{v_i}^*(v_i).
		\end{align*}
		Therefore, $\xi$ is one-harmonic if and only if $\operatorname{tr}(\mathcal{L}_\xi \nabla) = 0$, which gives the stated formula.
	\end{proof}
	\vspace{0.5cm}

	\begin{proof}[Proof of Theorem \ref{thm:main2}]
		If $\xi$ is a Killing vector field, then by Theorem \ref{thm:main1}, $\xi$ lies in the center of $\mathfrak{g}$, so $\operatorname{ad}_\xi = 0$. This immediately implies that equation \eqref{eq:one-harmonic} is satisfied, so $\xi$ is one-harmonic.
		
		Conversely, assume that $\xi$ is one-harmonic, and let $\{v_1, \ldots, v_n\}$ be an orthonormal basis of $\mathfrak{g}$. By Lemma \ref{lem:one-harmonic-char}, equation \eqref{eq:one-harmonic} holds. We take the inner product of both sides with $\xi$.
		
		\textbf{For the $\operatorname{ad}_{v_i}^*$-term:}
		\[
		\langle \operatorname{ad}_{v_i}^*(\operatorname{ad}_\xi v_i), \xi \rangle = \langle \operatorname{ad}_\xi v_i, \operatorname{ad}_{v_i} \xi \rangle = \langle \operatorname{ad}_\xi v_i, [v_i, \xi] \rangle = -\|\operatorname{ad}_\xi v_i\|^2,
		\]
		since $[v_i, \xi] = -\operatorname{ad}_\xi v_i$. Hence
		\[
		\sum_{i=1}^n \langle \operatorname{ad}_{v_i}^*(\operatorname{ad}_\xi v_i), \xi \rangle = -\sum_{i=1}^n \|\operatorname{ad}_\xi v_i\|^2.
		\]
		
		\textbf{For the $\operatorname{J}_{v_i}$-term:}
		\[
		\langle \operatorname{J}_{v_i}(\operatorname{ad}_\xi v_i), \xi \rangle = \langle \operatorname{ad}_{\operatorname{ad}_\xi v_i}^* v_i, \xi \rangle = \langle v_i, \operatorname{ad}_{\operatorname{ad}_\xi v_i} \xi \rangle = \langle v_i, \operatorname{ad}_\xi^2 v_i \rangle,
		\]
		where in the last equality, we use that $$\operatorname{ad}_{\operatorname{ad}_\xi v_i} \xi = [[\xi, v_i], \xi] = -[\xi, [\xi, v_i]] = -\operatorname{ad}_\xi^2 v_i.$$
		Summing over $i$ yields
		\[
		\sum_{i=1}^n \langle \operatorname{J}_{v_i}(\operatorname{ad}_\xi v_i), \xi \rangle = \operatorname{tr}(\operatorname{ad}_\xi^2).
		\]
		Since $\mathfrak{g}$ is nilpotent, the endomorphism $\operatorname{ad}_\xi$ is nilpotent, and hence $\operatorname{tr}(\operatorname{ad}_\xi^2) = 0$.
		
		\textbf{For the final term:} Set
		\[
		w := \sum_{i=1}^n \operatorname{ad}_{v_i}^*(v_i) \in \mathfrak{g}.
		\]
		We show that $w = 0$. For any $u \in \mathfrak{g}$,
		\[
		\langle w, u \rangle = \sum_{i=1}^n \langle \operatorname{ad}_{v_i}^*(v_i), u \rangle = \sum_{i=1}^n \langle v_i, \operatorname{ad}_{v_i} u \rangle = -\operatorname{tr}(\operatorname{ad}_u).
		\]
		Since $\mathfrak{g}$ is nilpotent, $\operatorname{ad}_u$ is nilpotent for all $u$, and therefore $\operatorname{tr}(\operatorname{ad}_u) = 0$. Thus $\langle w, u \rangle = 0$ for all $u$, which implies $w = 0$, and hence $\operatorname{ad}_\xi(w) = 0$.
		
		Combining these computations, taking the inner product of \eqref{eq:one-harmonic} with $\xi$ gives:
		\[
		0 = -\sum_{i=1}^n \|\operatorname{ad}_\xi v_i\|^2 + \operatorname{tr}(\operatorname{ad}_\xi^2) = -\sum_{i=1}^n \|\operatorname{ad}_\xi v_i\|^2.
		\]
		Therefore $\operatorname{ad}_\xi v_i = 0$ for every $i$, and hence $\operatorname{ad}_\xi = 0$. Thus $\xi$ lies in the center of $\mathfrak{g}$, and by Theorem \ref{thm:main1}, $\xi$ is a Killing vector field.
	\end{proof}
	
	\subsection{Conformal vector fields}
	An algebraic description of conformal vector fields on $G$
	is developed in the next lemma, which plays a foundational role in establishing Theorem \ref{thm:main3}.
	
	\begin{Lemma}\label{lem:conformal-char}
		A left-invariant vector field $\xi$ on $(G,h)$ is conformal if and only if 
		\[
		\mathcal{L}_{\xi}h=-\tfrac{2}{n}\operatorname{tr}\big(\operatorname{ad}_{\xi}\big)\langle\cdot,\cdot\rangle.
		\]
	\end{Lemma}
	
	\begin{proof}
		Let $\xi$ be a left-invariant vector field on $(G,h)$. Denote $\xi=\sum_{i=1}^n\xi_iv_i$ in an orthonormal basis $\{v_1,\cdots,v_n\}$ of $\mathfrak{g}$.
		
		\begin{align*}
			\mathrm{div}(\xi)&=\sum_{i=1}^n\langle \nabla_{v_i}\xi,v_i\rangle \\ 
			&=\sum_{i=1}^n\langle \mathrm{R}_{\xi}v_i,v_i\rangle\\
			&=-\tfrac{1}{2} \sum_{i=1}^n\langle( \operatorname{ad}_\xi + \operatorname{ad}_\xi^*  + \operatorname{J}_\xi)(v_i),v_i\rangle\\
			&=-\tfrac{1}{2} \sum_{i=1}^n\left[\langle\operatorname{ad}_\xi(v_i),v_i\rangle+\langle \operatorname{ad}_\xi^*(v_i),v_i\rangle+\langle \operatorname{J}_\xi(v_i),v_i\rangle\right]\\
			&=-\tfrac{1}{2}\Big(\operatorname{tr}(\operatorname{ad}_\xi)+\operatorname{tr}(\operatorname{ad}_\xi^*)+\operatorname{tr}(\operatorname{J}_\xi)\Big).
		\end{align*}
		Since $\operatorname{J}_\xi$ is skew-symmetric, $\operatorname{tr}(\operatorname{J}_\xi)=0$. Also, $\operatorname{tr}(\operatorname{ad}_\xi^*) = \operatorname{tr}(\operatorname{ad}_\xi)$.
		
		Therefore, $\mathrm{div}(\xi)=-\operatorname{tr}(\operatorname{ad}_\xi)$. So $\xi$ is conformal if and only if 
		\[
		\mathcal{L}_{\xi}h=\frac{2\,\mathrm{div}(\xi)}{n}\langle\cdot,\cdot\rangle = -\tfrac{2}{n}\operatorname{tr}\big(\operatorname{ad}_{\xi}\big)\langle\cdot,\cdot\rangle.
		\]
	\end{proof}

	\begin{proof}[Proof of Theorem \ref{thm:main3}]
		Let $\xi$ be a left-invariant conformal vector field on $(G,h)$. By Lemma \ref{lem:conformal-char},
		\[
		\mathcal{L}_{\xi}h=-\tfrac{2}{n}\operatorname{tr}\big(\operatorname{ad}_{\xi}\big)\langle\cdot,\cdot\rangle.\]
		Since $\mathfrak{g}$ is nilpotent, $\operatorname{ad}_\xi$ is nilpotent for any $\xi \in \mathfrak{g}$, and therefore from Lemma \ref{lem:nilpotent-skew}, $\operatorname{tr}(\operatorname{ad}_\xi)=0$.
		Thus, $\mathcal{L}_{\xi}h=0$, which means $\xi$ is a Killing vector field.\\
		
		Inversely, if $\xi$ is a Killing vector field on $G$, then $\operatorname{ad}_\xi=0$. Hence, $\operatorname{tr}\big(\operatorname{ad}_{\xi}\big)=0.$ It follows that $\mathcal{L}_{\xi}h=-\tfrac{2}{n}\operatorname{tr}\big(\operatorname{ad}_{\xi}\big)\langle\cdot,\cdot\rangle=0$. Accordingly $\xi$ is conformal.
	\end{proof}
	
	\begin{Remark}
		Since nilpotent Lie groups are unimodular, Theorem \ref{thm:main3} is a special case of the result in \cite{Cintra2016}.
	\end{Remark}
	
	\subsection{Concurrent vector fields}
	
	To prepare the ground for Theorem \ref{thm:main4}, the next lemma establishes the algebraic structure underlying concurrent vector fields on $G$.
	\begin{Lemma}\label{lem:concurrent-char}
		A left-invariant vector field $\xi$ on $(G,h)$ is concurrent if and only if $\mathrm{R}_{\xi}=\mathrm{id}_{\mathfrak{g}}$.
	\end{Lemma}
	
	\begin{proof}
		Let $\xi$ be a left-invariant vector field on $G$. Then $\xi$ is concurrent if and only if $\nabla_v(\xi)=v$ for all $v\in \mathfrak{g}$. This is equivalent to $\mathrm{R}_{\xi}(v)=v$ for all $v\in\mathfrak{g}$, i.e., $\mathrm{R}_\xi=\mathrm{id}_\mathfrak{g}$.
	\end{proof}

	\begin{proof}[Proof of Theorem \ref{thm:main4}]
		Suppose that a nontrivial left-invariant vector field $\xi$ on $(G,h)$ is concurrent. By Lemma \ref{lem:concurrent-char}, $\mathrm{R}_\xi = \mathrm{id}_\mathfrak{g}$.
		
		Recall that 
		\[
		\mathrm{R}_\xi = -\tfrac{1}{2} \left( \operatorname{ad}_\xi + \operatorname{ad}_\xi^* \right) - \tfrac{1}{2} \operatorname{J}_\xi.
		\]
		
		In particular, taking $v = \xi$, we have $\mathrm{R}_\xi(\xi) = \xi$, which gives:
		\[
		-\tfrac{1}{2} \left( \operatorname{ad}_\xi(\xi) + \operatorname{ad}_\xi^*(\xi) \right) - \tfrac{1}{2} \operatorname{J}_\xi(\xi) = \xi.
		\]
		
		Since $\operatorname{ad}_\xi(\xi)= 0$ and $\operatorname{J}_\xi(\xi) = \operatorname{ad}_\xi^*(\xi)$, the equation simplifies to:
		\[
		-\tfrac{1}{2} \operatorname{ad}_\xi^*(\xi) - \tfrac{1}{2} \operatorname{ad}_\xi^*(\xi) = \xi,
		\]
		which gives:
		\[
		-\operatorname{ad}_\xi^*(\xi) = \xi.
		\]
		
		Now, taking the inner product of both sides with $\xi$:
		\[
		-\langle \operatorname{ad}_\xi^*(\xi), \xi \rangle = \langle \xi, \xi \rangle = \|\xi\|^2.
		\]
		
		On the other hand, by definition of the adjoint operator:
		\[
		\langle \operatorname{ad}_\xi^*(\xi), \xi \rangle = \langle \xi, \operatorname{ad}_\xi(\xi) \rangle = 0.
		\]
		
		Therefore:
		\[
		\|\xi\|^2 = -\langle \operatorname{ad}_\xi^*(\xi), \xi \rangle = 0,
		\]
		which implies $\xi = 0$. This contradicts the assumption that $\xi$ is nontrivial.
		
		Hence, there exists no nontrivial concurrent left-invariant vector field on $(G,h)$.
	\end{proof}
	
	\begin{Remark}
		The proof shows that the concurrent condition $\nabla_v \xi = v$ for all $v$ leads to an algebraic constraint that cannot be satisfied by any nonzero vector in the Lie algebra. This is independent of the nilpotency assumption and holds for any Lie group with a left-invariant metric.
	\end{Remark}
	
	\begin{Corollary}
		On any Lie group $(G,h)$ with left-invariant Riemannian metric (not necessarily nilpotent), there are no nontrivial left-invariant concurrent vector fields.
	\end{Corollary}

	\subsection{Harmonic vector fields}
	The characterization rests on Lemma~\ref{prop:ricci-zero-central}, which establishes the vanishing of $\operatorname{Ric}(\xi)$ under the two conditions of Theorem \ref{thm:harmonic-rigorous}.
	\begin{Lemma}\label{prop:ricci-zero-central}
		Let $(G,h)$ be a connected nilpotent Lie group equipped with a left-invariant Riemannian metric, let $\xi$ be in the center of $\mathfrak{g}$ and  $\xi \perp [\mathfrak{g}, \mathfrak{g}]$. Then,
		\[
		\operatorname{Ric}(\xi) = 0.
		\]
	\end{Lemma}
	
	\begin{proof}
		Let $\xi$ be in the center of $\mathfrak{g}$ and  $\xi \perp [\mathfrak{g}, \mathfrak{g}]$.
		
		By Corollary 29  from \cite{Fok}, for a nilpotent Lie group, the Ricci curvature tensor is given by:
		\begin{equation}\label{Ric}
			\operatorname{ric}(u, v) = -\tfrac{1}{2} \operatorname{tr}(\operatorname{ad}_u \circ \operatorname{ad}_v^*) - \tfrac{1}{4} \operatorname{tr}(J_u \circ J_v).
		\end{equation}
		
		The Ricci operator $\operatorname{Ric}$ is defined by:
		\[
		\operatorname{ric}(u, v) = \langle \operatorname{Ric}(u), v \rangle.
		\]
		
		To show $\operatorname{Ric}(\xi) = 0$, we need to prove that $\langle \operatorname{Ric}(\xi), v \rangle = 0$ for all $v \in \mathfrak{g}$.
		
		This is equivalent to showing:
		\[
		\operatorname{ric}(\xi, v) = 0 \quad \forall v \in \mathfrak{g}.
		\]
		
		
		
		\textbf{Let's show that $\operatorname{tr}(\operatorname{ad}_\xi \circ \operatorname{ad}_v^*) = 0$.}
		
		Since $\xi$ is central:
		\[
		\operatorname{ad}_\xi(w) =  0 \quad \forall w \in \mathfrak{g}.
		\]
		
		Therefore, $\operatorname{ad}_\xi = 0$ as an endomorphism of $\mathfrak{g}$.
		
		Hence:
		\[
		(\operatorname{ad}_\xi \circ \operatorname{ad}_v^*)(w) = \operatorname{ad}_\xi(\operatorname{ad}_v^*(w)) = 0 \quad \forall w \in \mathfrak{g}.
		\]
		
		This means $\operatorname{ad}_\xi \circ \operatorname{ad}_v^* = 0$ as an endomorphism.
		
		Therefore:
		\[
		\operatorname{tr}(\operatorname{ad}_\xi \circ \operatorname{ad}_v^*) = 0.
		\]
		
		\textbf{ We now prove that $\operatorname{tr}(J_\xi \circ J_v) = 0$.}
		
		Let $\{v_1, \ldots, v_n\}$ be an orthonormal basis of $\mathfrak{g}$.
		
		By definition:
		\[
		\operatorname{tr}(J_\xi \circ J_v) = \sum_{i=1}^n \langle (J_\xi \circ J_v)(v_i), v_i \rangle.
		\]
		
		Using the trace property for endomorphisms with respect to an inner product:
		\[
		\operatorname{tr}(J_\xi \circ J_v) = \sum_{i=1}^n \langle J_v(v_i), J_\xi^*(v_i) \rangle.
		\]
		
		Now, we need to compute $J_\xi^*(v_i)$.
		
		
		
		For any $i,j\in \{1,\cdots,n\}$.
		
		
		
		Since $\xi\perp [\mathfrak{g},\mathfrak{g}]$,
		\[
		\langle J_\xi^*(v_i),v_j \rangle=\langle J_\xi(v_j) , v_i \rangle = \langle \operatorname{ad}_{v_j}^*(\xi), v_i \rangle = \langle \xi, \operatorname{ad}_{v_j}(v_i) \rangle=0.
		\]
		
		
		
		
		Then, 
		
		\[
		J_\xi^*(v_i) = 0.
		\]
		
		Therefore:
		\[
		\operatorname{tr}(J_\xi \circ J_v) = \sum_{i=1}^n \langle J_v(v_i), 0 \rangle = 0.
		\]
		
		Combining the two previous steps, we
		conclude:
		\[
		\operatorname{Ric}(\xi) = 0.
		\]
	\end{proof}

	\begin{proof}[Proof of Theorem \ref{thm:harmonic-rigorous}]
		
		\textbf{Step 1: (i) and (ii) imply harmonic.}
		
		Assume that $\xi$ is in the center of $\mathfrak{g}$ and  $\xi \perp [\mathfrak{g}, \mathfrak{g}]$.
		By Theorem \ref{thm:main2}, any left-invariant Killing vector field on a nilpotent Lie group is one-harmonic, i.e.
		\[
		\Delta \xi = \operatorname{Ric}(\xi).
		\]

		And moreover from Lemma \ref{prop:ricci-zero-central} a central vector field $\xi$ satisfying (ii), one has
		\[
		\operatorname{Ric}(\xi) = 0.
		\]

		Therefore, for a left-invariant field $\xi$ satisfying (i) and (ii), we have

		\[
		\Delta \xi = \operatorname{Ric}(\xi) = 0,
		\]

		So $\xi$ is harmonic.
		
		\medskip
		
		\textbf{Step 2: Harmonic implies (i) and (ii).}
		
		Now assume that $\xi$ is harmonic, i.e.\ $\Delta \xi = 0$.
		Let $\{v_1,\dots,v_n\}$ be an orthonormal basis of $\mathfrak{g}$.
		The Hodge--de Rham Laplacian on vector fields is given by

		\[
		\Delta \xi
		= -\sum_{i=1}^n \big( \operatorname{L}_{v_i}\operatorname{L}_{v_i}\xi - \operatorname{L}_{\operatorname{L}_{v_i}v_i}\xi \big).
		\]

		Taking the inner product with $\xi$ and using $\Delta \xi = 0$, we obtain

		\[
		0 = \langle \Delta \xi, \xi \rangle
		= -\sum_{i=1}^n \langle \operatorname{L}_{v_i}\operatorname{L}_{v_i}\xi, \xi \rangle
		+ \sum_{i=1}^n \langle \operatorname{L}_{\operatorname{L}_{v_i}v_i}\xi, \xi \rangle.
		\]

		For each $i$, by the Leibniz rule and metric compatibility,

		\[
		\langle \operatorname{L}_{v_i}\operatorname{L}_{v_i}\xi, \xi \rangle
		= v_i \langle \operatorname{L}_{v_i}\xi, \xi \rangle - \langle \operatorname{L}_{v_i}\xi, \operatorname{L}_{v_i}\xi \rangle.
		\]

		Moreover,

		\[
		\langle \operatorname{L}_{v_i}\xi, \xi \rangle
		= \tfrac{1}{2} v_i \langle \xi,\xi \rangle.
		\]

		Since $\xi$ is left-invariant, its norm $\|\xi\|^2 = \langle \xi,\xi \rangle$ is constant, hence

		\[
		v_i \langle \xi,\xi \rangle = 0
		\quad\text{and}\quad
		v_i \langle \operatorname{L}_{v_i}\xi, \xi \rangle = 0.
		\]

		Thus

		\[
		\langle \operatorname{L}_{v_i}\operatorname{L}_{v_i}\xi, \xi \rangle
		= -\|\operatorname{L}_{v_i}\xi\|^2.
		\]

		For the second term, again by metric compatibility,

		\[
		\langle \operatorname{L}_{\operatorname{L}_{v_i}v_i}\xi, \xi \rangle
		= \tfrac{1}{2} (\operatorname{L}_{v_i}v_i)\langle \xi,\xi \rangle = 0,
		\]

		since $\|\xi\|^2$ is constant.
		
		Putting everything together, we get

		\[
		0 = \langle \Delta \xi, \xi \rangle
		= -\sum_{i=1}^n \|\operatorname{L}_{v_i}\xi\|^2.
		\]
		
		Hence
		\begin{equation}\label{eqep}
			\operatorname{L}_{v_i}\xi = 0 \quad \forall i,
		\end{equation}

		and by linearity,

		\[
		\operatorname{L}_v \xi = 0 \quad \text{for all } v \in \mathfrak{g}.
		\]

		Thus $\xi$ is a parallel left-invariant vector field.
		
		\medskip
		
		
		For all $v,w \in \mathfrak{g}$, the Koszul formula for left-invariant vector fields gives

		\[
		2\langle \operatorname{L}_v \xi, w \rangle
		= \langle [v,\xi], w \rangle - \langle [\xi,w], v \rangle + \langle [w,v], \xi \rangle.
		\]

		Since $\operatorname{L}_v \xi = 0$ for all $v$, we obtain

		\[
		0
		= \langle [v,\xi], w \rangle - \langle [\xi,w], v \rangle + \langle [w,v], \xi \rangle
		\quad \text{for all } v,w \in \mathfrak{g}.
		\]

		Using $[\xi,w] = -[w,\xi]$, this can be rewritten as

		\[
		\langle [v,\xi], w \rangle + \langle [w,\xi], v \rangle
		= \langle [v,w], \xi \rangle.
		\]

		The right-hand side is antisymmetric in $(v,w)$, while the left-hand side is symmetric in $(v,w)$.
		Exchanging $v$ and $w$ and adding the two equations, we get

		\[
		2\big( \langle [v,\xi], w \rangle + \langle [w,\xi], v \rangle \big) = 0,
		\]

		so

		\[
		\langle [v,\xi], w \rangle + \langle [w,\xi], v \rangle = 0
		\quad \text{for all } v,w \in \mathfrak{g}.
		\]

		This means that the endomorphism

		\[
		\operatorname{ad}_\xi : \mathfrak{g} \to \mathfrak{g}, \quad v \mapsto [\xi,v],
		\]

		is skew-adjoint with respect to the inner $\langle \cdot, \cdot\rangle$.
		
		On the other hand, since $\mathfrak{g}$ is nilpotent, $\operatorname{ad}_\xi$ is a nilpotent endomorphism.
		
		Then by Lemma \ref{lem:nilpotent-skew}

		\[
		\operatorname{ad}_\xi = 0,
		\]

		Therefore, $\xi$ is central.
		
		
		
		It follows from the Koszul formula that,
		\[
		2\langle \operatorname{L}_{v_i}\xi, v_j \rangle = 0 - \langle [v_i, v_j], \xi \rangle - 0 = -\langle [v_i, v_j], \xi \rangle.
		\]
		
		Since $\operatorname{L}_{v_i}\xi = 0$ (equation \eqref{eqep}), we have $\langle \operatorname{L}_{v_i}\xi, v_j \rangle = 0$ for all $j$.
		
		Therefore:
		\[
		\langle [v_i, v_j], \xi \rangle = 0 \quad \forall i, j.
		\]
		
		Since $[\mathfrak{g}, \mathfrak{g}] = \mathrm{span}\{[v_i, v_j] : i, j = 1, \ldots, n\}$, this means:
		\[
		\xi \perp [\mathfrak{g}, \mathfrak{g}].
		\]
		
		This completes the proof.

		
		
		
	\end{proof}

	\end{document}